\documentclass{amsart}
\usepackage{amsmath,latexsym}
\usepackage[mathscr]{eucal}
\usepackage{amssymb}

\numberwithin{equation}{section}
 \allowdisplaybreaks \theoremstyle{remark}

\def\star{\raise .5ex \hbox{*}}
\def\sumstar_#1{\setbox0=\hbox{$\scriptstyle{#1}$}
\setbox2=\hbox{$\displaystyle{\sum}$}
\setbox4=\hbox{${}\star\mathsurround=0pt$} \dimen0=.5\wd0
\advance\dimen0 by-.5\wd2 \ifdim\dimen0>0pt \ifdim\dimen0>\wd4
\kern\wd4 \else\kern\dimen0\fi\fi \mathop{{\sum}\star}_{\kern-\wd4
#1}}

\begin{document}

\title[Primes in Tuples IV]{Primes in Tuples IV: Density of small gaps between
consecutive primes}
\author{D. A. Goldston, J. Pintz and C. Y. Y{\i}ld{\i}r{\i}m }

\thanks{The first author was supported by NSF; this study was done while
the second and third authors were members at the Institute for
Advanced Studies, Princeton during Fall 2009 and were supported by
the Oswald Veblen Fund; the second author also acknowledges the
partial support of ERC-AdG. No.228005.}

\date{\today}

\maketitle

\vskip .3in
\section{Introduction }
Let, as usual, $\pi(x)$ denote the number of primes $\leq x$, and
$p_{n}$ the $n$-th prime. The prime number theorem says $\pi(x)\sim
{x\over \log x}$ as $x\to\infty$, so that on average $p_{n+1}-p_{n}$
is $\log p_n$. In \cite{GPY1} we proved that
\begin{equation}
\liminf_{n\to\infty}{p_{n+1}-p_{n}\over \log p_{n}} =0. \label{eq:
1.1}
\end{equation}

The main result of this paper is \vskip .1in

\noindent {\bf Theorem 1} \,\, {\it For any fixed $\eta >0$, we
have}
\begin{equation}
\#\{p_n \leq x \, ; \, p_{n+1}-p_n \leq \eta\log p_n \} \gg_{\eta}
\pi(x), \label{eq: 1.2}
\end{equation}
{\it i.e. the small gaps between primes attested by the proof of
(1.1) in fact constitute a positive proportion of the set of all
gaps between consecutive primes}.  \vskip .1in

The course of the proof leads to Theorem 2 involving the explicit
estimate (3.26) below. In the last section we conditionally find
stronger versions of (3.26), based upon assuming more than what the
Bombieri-Vinogradov theorem provides. These results also provide
within the framework of our method a new quantitative manifestation
of the effect of the extent of assumed information on how well the
primes are distributed in arithmetic progressions. Theorem 3 in the
fourth section expresses that very small gaps are sparse. A rather
qualitative concise version of the results of this paper has been
presented in \cite{GPY}.

Concerning our subject matter it is conjectured (see \cite{So1},
\cite{So2}) that, given $0 \leq \alpha < \beta$, as $x\to\infty$ we
have
\begin{equation}
\#\{p_n \leq x \, ; \, p_{n+1}\in (p_n +\alpha\log p_n,p_n
+\beta\log p_n) \} \sim \pi(x)\int_{\alpha}^{\beta}e^{-t}\, dt.
\label{eq: 1.3}
\end{equation}
Gallagher's \cite{Ga} calculation shows that this conjecture can be
deduced from the Hardy-Littlewood prime $k$-tuples conjecture.

A bibliography of former results for the limit in (1.1) was given in
\cite{GPY1}. Before \cite{GPY1} the best known result for this limit
was
\begin{equation}
\liminf_{n\to\infty}{p_{n+1}-p_{n}\over \log p_{n}} \leq
0.2484\ldots \, , \label{eq: 1.4}
\end{equation}
due to Maier \cite{M}. In his proof Maier employed specially
constructed thin sets in which the density of primes is larger by a
factor of $e^{\gamma}$ than on average, and paralleling the
Bombieri-Davenport-Huxley method (\cite{BD}, \cite{H}) with
necessary modifications, he attained (1.4) which is $e^{-\gamma}$
times the result of this former method. Since Maier's method
involved working within a thin set of integers, the small gaps
indicated by (1.4) could not provide a positive proportion of all
gaps. The second best result was $\leq {1\over 4}$ by Goldston and
Y{\i}ld{\i}r{\i}m \cite{GY}. The origin of its method also being the
Bombieri-Davenport proof, the small gaps found in \cite{GY} were
shown to occur in a positive proportion of all cases.

Let us recall briefly how the result (1.1) was obtained. Consider
the $k$-tuple
\begin{equation}
\mathcal H = \{h_{1}, h_{2}, \ldots, h_{k}\} \;\; {\rm with \,\,
distinct \,\, integers} \;\; h_{1}, \ldots, h_{k} \in [1,h]  ,
\label{eq: 1.5}
\end{equation}
and for a prime $p$ denote by $\nu_{p}(\mathcal H)$ the number of
distinct residue classes modulo $p$ occupied by the entries of
$\mathcal H$. The singular series associated with $\mathcal H$ is
defined as
\begin{equation}
\mathfrak S (\mathcal H) : = \prod_{p}(1-{1\over
p})^{-k}(1-{\nu_{p}(\mathcal H)\over p}) , \label{eq: 1.6}
\end{equation}
the product being convergent because $\nu_{p}(\mathcal H)=k$ for
$p>h$. We say that $\mathcal H$ is admissible if
\begin{equation}
P_{\mathcal H}(n) := (n+h_{1})(n+h_{2})\cdots (n+h_{k}) \label{eq:
1.7}
\end{equation}
is not divisible by a fixed prime number for every $n$, which is
equivalent to $\nu_{p}(\mathcal H) \neq p$ for all $p$ and therefore
also to $\mathfrak S(\mathcal H) \neq 0$. That $\{ n+h_1 , n+h_2 ,
\ldots, n+h_k \}$ is a prime tuple, i.e. each entry is prime, is
equivalent to $P_{\mathcal H}(n)$ being a product of $k$ primes.
Since the generalized von Mangoldt function
\begin{displaymath}
\Lambda_{k}(m) := \sum_{d \mid m}\mu(d)(\log {m\over d})^k
\end{displaymath}
vanishes when $m$ has more than $k$ distinct prime factors, the
quantity
\begin{displaymath}
{1\over k!}\sum_{\substack{d\mid P_{\mathcal H}(n) \\ d \leq
R}}\mu(d)(\log {R\over d})^k
\end{displaymath}
with the truncation $d\leq R$ may be employed in a roughly correct
detection of prime tuples (the contribution from proper prime power
factors is negligible; $1/k!$ is a normalization factor). One
crucial idea is to give up trying to count tuples consisting of
primes exclusively, but rather also include tuples with primes in
many entries. This brings about the use of
\begin{equation}
\Lambda_{R}(n; \mathcal H, \ell) := {1\over
(k+\ell)!}\sum_{\substack{d\mid P_{\mathcal H}(n) \\ d \leq
R}}\mu(d)(\log {R\over d})^{k+\ell} , \quad (0 \leq \ell < k)
\label{eq: 1.8}
\end{equation}
so as to allow counting those $P_{\mathcal H}(n)$ which have at most
$k+\ell$ distinct prime factors.

Let
\begin{equation}
\theta(n) := \begin{cases}
\log n &\text{ if $n$ is prime}, \\
0 &\text{ otherwise},
\end{cases}
\label{eq: 1.9}
\end{equation}
and
\begin{equation}
\Theta(n,h) := \sum_{1 \leq h_{0} \leq h}\theta(n+h_{0}) .
\label{eq: 1.10}
\end{equation}
The proof of (1.1) is achieved by showing the positivity of the
quantity
\begin{equation}
S_{R}(N,k,\ell,h) := \sum_{N< n \leq 2N}\left(\Theta(n,h) - \log
3N\right) \bigl(\sum_{\scriptstyle \mathcal H \subset [1, h] \atop
\scriptstyle |\mathcal H | = k} \Lambda_{R}(n;\mathcal H ,
\ell)\;\bigr)^{2}. \label{eq: 1.11}
\end{equation}
Here, as $N\to\infty$, for a result of the type (1.1) we need
$\epsilon\log N \ll h \ll \log N$ with an arbitrarily small but
fixed $\epsilon > 0$, and the larger the truncation level $R$ is
relative to $N$ the better detection will be provided by (1.8). The
tuple size $k$ is taken to be arbitrarily large but fixed. In fact
for the proof of (1.1) it suffices to consider the simpler
expression where the inner sum consists only of the diagonal terms
$\Lambda_{R}^{2}(n;\mathcal H , \ell)$, and a modified version of
this will be used in Section 3. The expression in (1.11) is needed
for achieving a better result in the case of the gaps $p_{n+r}-p_n$
with $r \geq 2$ in \cite{GPY1} and for a stronger quantitative
version of (1.1) in \cite{GPY2}.

The information on primes, beyond the prime number theorem, that is
of key importance in our studies is the level of distribution of
primes in arithmetic progressions. We say that the primes satisfy a
level of distribution $\vartheta$ if
\begin{equation}
\sum_{q\leq Q}\max_{\scriptstyle a \atop \scriptstyle (a,q)=1}\Bigl|
\sum_{\scriptstyle p : {\rm prime}\ \atop {\scriptstyle p\leq N
\atop \scriptstyle p\equiv a (\!\bmod q)}}\log p - {N\over
\phi(q)}\Bigr| \ll_{\epsilon, A} {N\over (\log N)^A} \label{eq:
1.12}
\end{equation}
holds for any $A >0$ and any $\epsilon > 0$ with
\begin{equation}
Q = N^{\vartheta - \epsilon}. \label{eq: 1.13}
\end{equation}
According to the Bombieri-Vinogradov theorem, for any $A>0$ there is
a $B=B(A)$ such that (1.12) holds with $Q=N^{{1\over 2}}(\log
N)^{-B}$, so that the primes are known to have level of distribution
${1\over 2}$. The Elliott-Halberstam conjecture is that the primes
have level of distribution $1$.

The following are special cases of some results from \cite{GPY1}
which are relevant to our purpose in this article. For an admissible
$k$-tuple $\mathcal{H}$, we have
\begin{equation}
\sum_{n\le N} \Lambda_R( n;\mathcal{H},\ell)^2 \sim {2\ell \choose
\ell} \frac{(\log R)^{k + 2\ell}}{(k + 2\ell)!} \mathfrak S(\mathcal
H)N , \label{eq: 1.14}
\end{equation}
as $R, N\to \infty$, for $R\ll N^{\frac{1}{2}}(\log N)^{-8M}$ where
$M= k+\ell$, and $h\le R^C$ for any given constant $C > 0$. In the
situation of weighting with the primes, for $1\le h_0\le h$ writing
$m=1$ when $h_{0}\in \mathcal{H}$ and $m=0$ when $h_{0}\not\in
\mathcal{H}$, if $\mathcal H \cup \{h_{0}\}$ is admissible we have
\begin{equation}
\sum_{N \leq n} \theta (n+h_0)\Lambda_R (n;\mathcal{H},\ell)^2 \sim
{2(\ell +m) \choose \ell +m}\frac{\mathfrak S(\mathcal H \cup
\{h_{0}\})}{(k + 2\ell +m)!} N(\log R)^{k + 2\ell +m}. \label{eq:
1.15}
\end{equation}
as $R, N\to \infty$, provided that $R\ll_M N^{\frac{1}{4}}(\log
N)^{-B(M)}$ for a sufficiently large positive constant $B(M)$, and
$h\le R$. The upper bound for $R$ is forced by the dependence of the
proof of (1.15) on the Bombieri-Vinogradov theorem, and for the
unconditional results in \cite{GPY1}, taking $R=N^{{1\over
4}-\epsilon}$ suffices. More generally, (1.15) holds with $R\ll
N^{\frac{\vartheta}{2} - \epsilon}$ and $h\le R^\epsilon$ for any
$\epsilon>0$,  assuming that the primes have level of distribution
$\vartheta$ with a fixed $\vartheta \in [{1\over 2}, 1]$.

The proof of (1.14) and (1.15) may be outlined as follows. Upon
writing the left-hand sides explicitly by substituting (1.8), the
sum over $n$ is carried to the innermost position and easily
evaluated. Then a Mellin transform converts the expressions into
integrals over vertical lines in the complex plane. The integrands
contain Dirichlet series which encode the arithmetic information
from the tuples. The integrals are evaluated by shifting the lines
of integration appropriately and by calculating the residues and the
bounds for the integrals over the new contours.

For the calculation of $S_{R}(N,k,\ell,h)$, the general versions of
(1.14) and (1.15) are employed in the expression on the right-hand
side of (1.11), and then Gallagher's \cite{Ga} result
\begin{equation}
\sum_{\scriptstyle \mathcal H \subset  [1, h] \atop \scriptstyle
|\mathcal H|=k} \mathfrak S (\mathcal H) \sim h^k \qquad {\rm for
\;\, fixed} \;\, k \;\, {\rm as} \;\; h\to\infty \label{eq: 1.16}
\end{equation}
(where each set is counted $k!$ times due to all of its
permutations) is needed to complete the calculation. The parameters
which appear in this process are chosen judiciously, in particular
$k$ has to be arbitrarily large but fixed and the optimal order of
magnitude of the integer $\ell$ turns out to be $\sqrt{k}$.

The proof of the positive proportion result in \cite{GY} uses an
argument which depends on the calculation of the fourth moment of
prime tuple approximants. If that argument is adapted
straightforwardly to the approach which led to (1.1), for a proof of
Theorem 1 one needs to show that
\begin{equation}
\sum_{N <n \leq 2N} \Bigl(\sum_{\scriptstyle \mathcal H \subset [1,
h] \atop \scriptstyle |\mathcal H | = k} \Lambda_{R}(n;\mathcal H ,
\ell)\;\Bigr)^{4} \ll N(\log N)^{4k + 4\ell}. \label{eq: 1.17}
\end{equation}
However, upon some calculations, the truth of this seems to be
questionable.

\section{Some preliminaries}

The lack of success from a direct use of results from \cite{GPY1}
for the proof of a positive proportion result notwithstanding, a
version of (1.14) and (1.15) in which the $n$ with $P_{\mathcal
H}(n)$ having small prime factors are discounted vouchsafes the
solution. We define
\begin{equation}
\mathcal{P}(x) := \prod_{p_{n}\leq x}p_n . \label{eq: 2.1}
\end{equation}
We shall use the following results which are consequences of (1.14),
(1.15) and Lemmas 4 and 5 of Pintz's work \cite{P}. \vskip .1in

\noindent {\bf Proposition 1} \,\, {\it For $N^{c_1} \leq R \leq
N^{{1\over 2+\delta}}(\log N)^{-c_2}$ where $c_1$ and $c_2$ are
suitably chosen constants depending on $k$ and $\ell \asymp
\sqrt{k}$ ($c_1$ can be taken to be ${1\over 5}$ and $c_2$ is
sufficiently large), $\delta > 0$ small compared to $k^{-{3\over
2}}$, $\mathcal H$ admissible with $h \ll \log R$ and $h\to\infty$
with $N$, we have}
\begin{equation}
\sum_{\scriptstyle N < n\le 2N \atop \scriptstyle (P_{\mathcal
H}(n), \mathcal{P}(R^{\delta}))=1} \Lambda_R( n;\mathcal{H},\ell)^2
\sim \bigl(1+O(k^{3}\delta^{2})\bigr) {2\ell \choose \ell}
\frac{\mathfrak S(\mathcal H)}{(k + 2\ell)!} N(\log R)^{k + 2\ell}.
\label{eq: 2.2}
\end{equation}
\vskip .1in

\noindent {\bf Proposition 2} \,\, {\it Upon the conditions of
Proposition 1 and the notation introduced in connection with (1.15),
if the level of distribution of primes is $\vartheta \geq {1\over
2}$, then for $N^{c_1} \leq R \leq N^{{\vartheta - \epsilon\over
(2+\delta)}}(\log N)^{-c_2}$ ($\epsilon >0$ arbitrarily small but
fixed) and $\mathcal H \cup \{h_{0}\}$ admissible, we have}
\begin{align}
& \sum_{\scriptstyle N < n\le 2N \atop \scriptstyle (P_{\mathcal
H}(n), \mathcal{P}(R^{\delta}))=1} \theta
(n+h_0)\Lambda_R (n;\mathcal{H},\ell)^2 \nonumber \\
& \qquad \qquad \sim \bigl(1+O(k^{3}\delta^{2})\bigr){2(\ell +m)
\choose \ell +m}\frac{\mathfrak S(\mathcal H \cup \{h_{0}\})}{(k +
2\ell +m)!} N(\log R)^{k + 2\ell +m}; \label{eq: 2.3}
\end{align}
{\it in case $\mathcal H \cup \{h_{0}\}$ is not admissible, the
right-hand side of (2.3) is $o(N(\log R)^{k + 2\ell +m})$.}

\vskip .1in

\noindent {\it Proof:} \,\,\, From \cite{P}, along with (1.14),
(1.15) these results are obvious except that the present error term
$O(k^{3}\delta^{2})$ meant with an absolute constant comes out as
$O(\delta)$ with the constant implied depending on $k$ and $\ell$.
Since we shall use the dependence on $k$ and $\ell$ of the error
term, we give its proof. An examination of the proof of Pintz's
Lemma 3 reveals that we need to have more precise versions of
(6.17), (6.18), (6.25) and (6.26) of \cite{P}. First we evaluate
\begin{align}
{\mathcal T}_{q,1}(1+\alpha) & := {1\over \ell !}\Big[\bigl({d\over
d\xi}\bigr)^{\ell}\Bigl({(1+\alpha +\xi)^{k+2\ell}\over
(1+\xi)^{k}}\Bigr)\Big]_{\xi = 0} \nonumber \\
& = (1+\alpha)^{k+\ell}\sum_{m=0}^{\ell} {2\ell -m \choose
\ell}{k+m-1 \choose m} (-\alpha)^m . \label{eq: 2.4}
\end{align}
This indicates that it would be opportune to restrict $|\alpha|$ to
values small compared to ${1\over k}$. Assuming this, and recalling
that $\ell \asymp \sqrt{k}$, from (2.4) we see that
\begin{equation}
{\mathcal T}_{q,1}(1+\alpha) = {2\ell \choose \ell}\Bigl(1+({k\over
2}+\ell)\alpha + ({k^2\over 8}+{k\ell\over 2}-{3k\over
8}+{\ell^2\over 2}-{\ell\over 2}-{k(k+1)\over 8(2\ell
-1)}\bigr)\alpha^2 + O((k\alpha)^3)\Bigr). \label{eq: 2.5}
\end{equation}
We will denote the coefficient of $\alpha^2$ in the last line as
$K$. This is used in (6.8) of \cite{P}. We recall that the prime
number $q=R^{\beta}$ in the statement of Lemma 3. In the last factor
of the integrand of (6.8) there are four terms. With the notation
introduced in (6.11) of \cite{P} we have the following. The first
term has $R_1 = R_2 = R$, so that $\alpha = 0$, and we get the
contribution
\begin{displaymath}
{2\ell \choose \ell}{(\log R)^{k+2\ell}\over (k+2\ell)!}G_{q}(0,0).
\end{displaymath}
The second term has $R_1 = R/q,\, R_2 = R$, so that $\alpha =
-\beta$, and we get the contribution
\begin{displaymath}
{2\ell \choose \ell}{(\log R)^{k+2\ell}\over
(k+2\ell)!}G_{q}(0,0)\Big[1-({k\over 2}+\ell)\beta + K\beta^2 +
O((k\beta)^3)\Big].
\end{displaymath}
The third term has $R_1 = R,\, R_2 = R/q$, so that $\alpha =
{\beta\over 1-\beta}$, and we get the contribution
\begin{displaymath} {2\ell \choose \ell}{(\log
R^{1-\beta})^{k+2\ell}\over (k+2\ell)!}G_{q}(0,0)\Big[1+({k\over
2}+\ell){\beta\over 1-\beta} + K\bigl({\beta\over 1-\beta}\bigr)^2 +
O((k\beta)^3)\Big].
\end{displaymath} The fourth term has $R_1 = R_2 = R/q$, so that
$\alpha =0$, and we get the contribution \begin{displaymath} {2\ell
\choose \ell}{(\log R^{1-\beta})^{k+2\ell}\over
(k+2\ell)!}G_{q}(0,0). \end{displaymath} Combining these as in (6.8)
we obtain
\begin{equation}
{2\ell \choose \ell}{(\log R)^{k+2\ell}\over
(k+2\ell)!}G_{q}(0,0)\Big[\bigl({k^2\over 4}+k\ell +\ell^2 + {k\over
4}+{k(k+1)\over 4(2\ell -1)}\bigr)\beta^2 + O((k\beta)^3)\Big]
\label{eq: 2.6}
\end{equation}
in place of the main term of (6.25) of \cite{P}. Hence in the new
version for (6.1) of \cite{P}, instead of ${\beta\over q}$ we have
\begin{align}
& {\nu_{q}(\mathcal H)\over q}{G_{q}(0,0)\over
G(0,0)}\Big[\bigl({k^2\over 4}+k\ell +\ell^2 + {k\over
4}+{k(k+1)\over 4(2\ell -1)}\bigr)\beta^2
+ O((k\beta)^3)\Big] \nonumber \\
& \quad = {\nu_{q}(\mathcal H)\over q}\Bigl(1- {\nu_{q}(\mathcal
H)\over q}\Bigr)^{-1}\bigl({k^2 \beta^{2}\over 4} + {\rm smaller \;
terms}\bigr), \nonumber
\end{align}
i.e. we can express the new version of Lemma 3 of \cite{P} as
\begin{equation}
\sum_{\scriptstyle N < n\le 2N \atop \scriptstyle q | P_{\mathcal
H}(n)} \Lambda_R( n;\mathcal{H},\ell)^2 \leq {\nu_{q}(\mathcal
H)\over q-\nu_{q}(\mathcal H)}{k^2 \beta^{2}\over 3} \sum_{N < n\le
2N}\Lambda_R( n;\mathcal{H},\ell)^2 .\label{eq: 2.7}\end{equation}
Here $N^{c_1} \leq R \leq N^{{1\over 2+\delta}}(\log N)^{-c_2}$, $q$
is a prime number for which we write $q = R^{\beta}$, with $0< \beta
\leq \delta$ where $\delta$ is small compared to $k^{-{3\over 2}}$
say, $k$ is sufficiently large, and $\ell \asymp \sqrt{k}$.

We know that $\nu_{q}(\mathcal H) \leq \min(q-1,k)$ since $\mathcal
H$ is admissible. For $q \leq k$, we have $\nu_{q}(\mathcal H) \leq
q-1$, so that ${\nu_{q}(\mathcal H)\over q-\nu_{q}(\mathcal H)} \leq
q-1$. For $q
>k$, we take $\nu_{q}(\mathcal H) \leq k$, so that
${\nu_{q}(\mathcal H)\over q-\nu_{q}(\mathcal H)}\leq {k\over q-k}$.
Summing over all primes $q \leq R^{\delta}$ we obtain a new version
of Lemma 4 of \cite{P} as
\begin{equation}
\sum_{\scriptstyle N < n\le 2N \atop \scriptstyle (P_{\mathcal
H}(n), \mathcal{P}(R^{\delta}))> 1} \Lambda_R( n;\mathcal{H},\ell)^2
\leq {k^3 \delta^2\over 4}  \sum_{N < n\le 2N}\Lambda_R(
n;\mathcal{H},\ell)^2 \label{eq: 2.8} \end{equation} if $k$ is large
enough. We see that we have to choose $\delta$ small enough so that
$k^3 \delta^2$ will be small. Now by (1.14) and (1.15) we
immediately obtain (2.2). When there is the twisting with primes the
proof runs similarly and Proposition 2 also follows.

\section{Proof of Theorem 1; Theorem 2}
In order to prove Theorem 1 we need to show the inequality
\begin{equation}
\sum_{\scriptstyle N < p_{j} \leq 2N \atop \scriptstyle p_{j+1}-p_j
\leq h} 1 \gg \pi(N) \sim {N\over \log N}, \quad (N \to \infty),
\label{eq: 3.1}
\end{equation}
for
\begin{equation}
h=\eta\log N, \quad \!\! \eta >0 \, {\rm \; arbitrarily \; small \;
but \; fixed}.  \label{eq: 3.2}
\end{equation}

Let
\begin{equation}
Q(N,h):= \sum_{\scriptstyle N < n \leq 2N \atop \scriptstyle
\pi(n+h)-\pi(n)> 1} 1 . \label{eq: 3.3}
\end{equation}
If $n$ is an integer for which $\pi(n+h)-\pi(n) > 1$, then there
must be a $j$ such that $n < p_j$ and $p_{j+1} \leq n+h$. Thus
$p_{j+1}-p_j < h$ and $p_{j+1} -h \leq n < p_j$, so that there are
less than $\lfloor h \rfloor $ such integers $n$ corresponding to
each such gap. Therefore
\begin{equation}
Q(N,h) \leq \, h\sum_{\scriptstyle N < p_{j} \leq 2N \atop
\scriptstyle p_{j+1}-p_j \leq h}1 \, + O(Ne^{-c\sqrt{\log N}}),
\label{eq: 3.4}
\end{equation}
where we have used the prime number theorem with error term to
remove the prime gaps which overlap the endpoints. (This is
explicitly shown in \cite{GY1}).

Instead of $S_R$ which was defined in (1.11), we will work with
\begin{equation}
\tilde{S}_{R} := {1\over N(h\log R)^{k}}\sum_{N < n \leq
2N}\bigl(\Theta(n,h) - \log 3N\bigr) \bigl(\sumstar_{\mathcal H}
\Lambda_{R}^{2}(n;\mathcal H , \ell)\;\bigr), \label{eq: 3.5}
\end{equation}
where
\begin{equation}
\sumstar_{\mathcal H} := \sum_{\scriptstyle \mathcal H \subset  [1,
h],\, |\mathcal H|=k \atop {\scriptstyle \mathcal H : \,
\mathrm{admissible}\atop \scriptstyle (P_{\mathcal
H}(n),\mathcal{P}(R^{\delta}))=1}} \,. \label{eq: 3.6}
\end{equation}
We note that, as a function of $\eta$, $k$ and $\ell$ will be chosen
sufficiently large but fixed, and $\delta >0$ will be chosen
sufficiently small but fixed (see (3.22) below).

From (3.5) we have, when $N$ is sufficiently large,
\begin{align}
\tilde{S}_{R} & \leq {1\over N(h\log R)^{k}} \sum_{\scriptstyle N <
n \leq 2N \atop\scriptstyle \Theta(n,h) \geq {3\over 2}\log
N}\Theta(n,h) \sumstar_{\mathcal H}\Lambda_{R}^{2}(n;\mathcal H ,
\ell)  \label{eq: 3.7} \\
& \leq {1\over N(h\log R)^{k}} \{\sum_{\scriptstyle N < n \leq 2N
\atop \scriptstyle \Theta(n,h) \geq {3\over 2}\log N}1\}^{{1\over
2}} \{ \sum_{N <n\leq 2N}(\Theta(n,h))^2 \Bigl(\sumstar_{\mathcal
H}\Lambda_{R}^{2}(n;\mathcal H , \ell)\Bigr)^2 \}^{{1\over 2}} \nonumber \\
& = {Q(N,h)^{{1\over 2}}\over N(h\log R)^{k}} I^{{1\over 2}},
\nonumber
\end{align}
where
\begin{align}
I= \sum_{1\leq h', h'' \leq h}\sum_{\scriptstyle \mathcal{H}_i
\subset  [1, h],\, |\mathcal{H}_i|=k \atop {\scriptstyle
\mathcal{H}_i : \, \mathrm{admissible}\atop \scriptstyle i=1,2}}
\sum_{\scriptstyle N < n \leq 2N \atop \scriptstyle
(P_{\mathcal{H}_{1}\cup \mathcal{H}_{2}}(n),
\mathcal{P}(R^{\delta}))=1}& \theta(n+h')\theta(n+h'')\nonumber \\
& \; \times \Lambda_{R}^{2}(n;\mathcal{H}_{1},
\ell)\Lambda_{R}^{2}(n;\mathcal{H}_{2}, \ell).
\label{3.8}\end{align} Here for a number $n$ to make a nonzero
contribution both of $n+h'$ and $n+h''$ must be prime, so that
$\bigl((n+h')(n+h''), \mathcal{P}(R^{\delta})\bigr)=1$ and writing
$\mathcal{H}_0 = \{h'\} \cup \{h''\} \cup \mathcal{H}_{1}\cup
\mathcal{H}_{2}$ we can re-express the condition on $n$ as
$(P_{\mathcal{H}_{0}}(n), \mathcal{P}(R^{\delta}))=1$. We also
observe that, since all prime factors of $P_{\mathcal H}(n)$ in
$\displaystyle\sumstar_{\mathcal H}$ are greater than $R^{\delta}$,
the number of squarefree divisors of $P_{\mathcal H}(n)$ is at most
$2^{{k\log 3N\over \delta\log R}}$. So, for any term in
$\displaystyle\sumstar_{\mathcal H}$ we have
\begin{equation}
\Lambda_{R}(n;\mathcal H , \ell) \leq {2^{{k\log 3N\over \delta\log
R}}\over (k+\ell)!}(\log R)^{k+\ell}, \label{eq: 3.9}
\end{equation}
and therefore
\begin{equation}
I \leq {2^{{4k\log 3N\over \delta\log R}}(\log R)^{4(k+\ell)}(\log
3N)^2\over (k+\ell)!^4} \sum_{1\leq h', h'' \leq
h}\sum_{\scriptstyle \mathcal{H}_i \subset  [1, h],\,
|\mathcal{H}_i|=k \atop {\scriptstyle \mathcal{H}_i : \,
\mathrm{admissible}\atop \scriptstyle i=1,2}} \sum_{\scriptstyle N <
n \leq 2N \atop \scriptstyle (P_{\mathcal{H}_{0}}(n),
\mathcal{P}(R^{\delta}))=1}1 . \label{3.10}\end{equation} For a
given $\mathcal{H}_0 \subset \{1, \ldots , h\}$ with
$|\mathcal{H}_0| = k+r,\, 0 \leq r \leq k+2$, denoting by $D(k,r)$
the number of quadruples $h',h'',\mathcal{H}_1,\mathcal{H}_2$
corresponding to $\mathcal{H}_0$, we re-express (3.10) as
\begin{equation}
I  \leq {2^{{4k\log 3N\over \delta\log R}}(\log R)^{4(k+\ell)}(\log
3N)^2\over (k+\ell)!^4}\sum_{r=0}^{k+2}D(k,r)\sum_{|\mathcal{H}_0| =
k+r}\sum_{\scriptstyle N < n \leq 2N \atop \scriptstyle
(P_{\mathcal{H}_{0}}(n), \mathcal{P}(R^{\delta}))=1}1 .
\label{3.11}\end{equation} We now invoke the main theorem of
Selberg's upper bound sieve (Theorem 5.1 of \cite{HR} or Theorem 2
in \S 2.2.2 of \cite{Gr}) that for any set $\mathcal H$ and $\delta
< {1\over 2}$
\begin{equation}
\sum_{\scriptstyle N < n \leq 2N \atop \scriptstyle
(P_{\mathcal{H}}(n), \mathcal{P}(R^{\delta}))=1}1 \leq
{N|\mathcal{H}|!\mathfrak S (\mathcal{H})\over (\log
R^{\delta})^{|\mathcal{H}|}}(1+o(1)), \quad \; (N\to\infty),
\label{eq: 3.12}
\end{equation}
which gives upon using (1.16) that
\begin{align}
I & \lesssim N {2^{{4k\log 3N\over \delta\log R}}(\log
R)^{4(k+\ell)}(\log N)^2\over
(k+\ell)!^4}\sum_{r=0}^{k+2}{(k+r)!D(k,r)\over (\delta\log R)^{k+r}}
\sum_{\mathcal{H}_0;\, |\mathcal{H}_0| =
k+r}\mathfrak{S}(\mathcal{H}_0) \label{eq: 3.13} \\
& \lesssim N {2^{{4k\log 3N\over \delta\log R}}(\log
R)^{4(k+\ell)}(\log N)^2\over
(k+\ell)!^4}\sum_{r=0}^{k+2}(k+r)!D(k,r)\Bigl({h \over \delta\log
R}\Bigr)^{k+r}. \nonumber
\end{align}
To deal with the inner sum here, first note that
\begin{equation}
D(k,r) := {k!^{2}(k+r)! (k^4 + 3k^3 + (3r+2)k^2 + 4rk+r^2)\over r!^2
(k+2-r)!} \label{eq:3.14}
\end{equation}
(here the factor $k!^{2}$ comes from the ordering of the elements
within the $k$-tuples $\mathcal{H}_1$ and $\mathcal{H}_2$). We skip
the proof of (3.14) since it follows from an elementary
combinatorial calculation, and after our choice of parameters the
order of magnitude is much smaller than that of the heftiest factor
$ 2^{{4k\log 3N\over \delta\log R}}$. Now for $u>0$, we have
\begin{align}
\sum_{r=0}^{k+2}(k+r)!D(k,r)u^{k+r} & \leq u^k (k+1)^2 (k+2)^2
\sum_{r=0}^{k+2}{k!^2 (k+r)!^2\over r!^2 (k+2-r)!}u^r \nonumber \\
& = u^k (k+2)!\sum_{r=0}^{k+2}{(k+r)!^2\over r!} {k+2
\choose r}u^r \nonumber \\
& \leq (2k+2)!^2 u^k (1+u)^{k+2}\label{eq: 3.15} ,
\end{align}
so that
\begin{equation}
I \lesssim N (\log R)^{4(k+\ell)}(\log N)^2 {(2k+2)!^2 \over
(k+\ell)!^4} 2^{{4k\log 3N\over \delta\log R}}\Bigl({h \over
\delta\log R}\Bigr)^{k} \Bigl(1+ {h \over \delta\log R}\Bigr)^{k+2}.
\label{3.16}\end{equation}

Using (3.16) and (3.4) in (3.7), we obtain
\begin{align}
\tilde{S}_{R} \lesssim & \Big(h\sum_{\scriptstyle N < p_{j} \leq 2N
\atop \scriptstyle p_{j+1}-p_j \leq h}1 \, + O(Ne^{-c\sqrt{\log
N}})\Big)^{{1\over 2}} \nonumber\\
& \quad\times {(\log R)^{(k+2\ell)}\log N\over N^{{1\over
2}}h^{k}}{(2k+2)!\over (k+\ell)!^2} 2^{{2k\log 3N\over \delta\log
R}}\Bigl({h \over \delta\log R}\Bigr)^{{k\over 2}} \Bigl(1+ {h \over
\delta\log R}\Bigr)^{{k+2\over 2}}. \label{3.17}\end{align}

Now we calculate $\tilde{S}_{R}$ using Propositions 1 and 2. From
Proposition 1 and (1.16) we see that
\begin{align}
\sum_{N < n\le 2N}\log 3N  & \sumstar_{\mathcal H}\Lambda_R(
n;\mathcal{H},\ell)^2 \label{eq: 3.18} \\
& \sim \bigl(1+O(k^3 \delta^2)\bigr) {2\ell \choose \ell}
\frac{h^{k}}{(k + 2\ell)!} N(\log R)^{k + 2\ell}\log N . \nonumber
\end{align}
Similarly, Proposition 2 and (1.16) imply
\begin{align}
&\sum_{\scriptstyle \mathcal{H}\subset  [1, h], \, |\mathcal{H}|=k
\atop \scriptstyle \mathcal{H} : \, \mathrm{admissible} } \sum_{h_i
\in \mathcal{H}}\sum_{\scriptstyle N < n\le 2N \atop \scriptstyle
(P_{\mathcal H}(n), \mathcal{P}(R^{\delta}))=1} \theta
(n+h_i)\Lambda_R (n;\mathcal{H},\ell)^2 \nonumber \\
& \qquad \sim \bigl(1+O(k^3 \delta^2)\bigr){2\ell +2 \choose \ell
+1}\frac{kh^{k}}{(k + 2\ell +1)!} N(\log R)^{k + 2\ell +1},
\label{eq: 3.19}
\end{align}
and
\begin{align}
& \sum_{\scriptstyle \mathcal{H}\subset  [1, h], \, |\mathcal{H}|=k
\atop \scriptstyle \mathcal{H} : \, \mathrm{admissible}
}\sum_{\scriptstyle 1 \leq h_0 \leq h \atop h_0 \not\in\mathcal{H}}
\sum_{\scriptstyle N < n\le 2N \atop \scriptstyle (P_{\mathcal
H}(n), \mathcal{P}(R^{\delta}))=1} \theta
(n+h_0)\Lambda_R (n;\mathcal{H},\ell)^2 \nonumber \\
& \qquad \gtrsim \bigl(1+O(k^3 \delta^2)\bigr){2\ell \choose
\ell}\frac{h^{k+1}}{(k + 2\ell)!} N(\log R)^{k + 2\ell}. \label{eq:
3.20}
\end{align}
Putting (3.18)-(3.20) together in (3.5) we obtain
\begin{equation}
\tilde{S}_{R} \gtrsim {{2\ell \choose \ell}\over (k+2\ell)!}(\log N)
(\log R)^{2\ell}\Big\{{k\over k+2\ell+1}{2(2\ell +1)\over \ell
+1}{\log R\over \log N} + \eta -1 +O(k^3 \delta^2)\Big\}. \label{eq:
3.21}
\end{equation}
Now given a small fixed $\eta >0$ if we take
\begin{equation}
\ell = \lfloor{4\over \eta}\rfloor , \quad k=2(\ell +1)(2\ell +1),
\quad \delta = {1\over \ell^4}, \quad R= N^{{1\over 4(1+\delta)}},
\label{eq: 3.22}
\end{equation}
then for the factor in brackets in (3.21) we see that
\begin{equation}
\Big\{{k\over k+2\ell+1}{2(2\ell +1)\over \ell +1}{\log R\over \log
N} + \eta -1 +O(k^3 \delta^2)\Big\} > {\eta\over 2} \label{eq: 3.23}
\end{equation}
holds for sufficiently small $\eta$, so that $\tilde{S}_{R} >0$.

Noting that (3.2) and (3.22) imply ${h\over \delta\log R} =
4\eta(\ell^4 +1) > 16{\ell^4 +1\over \ell +1} \geq 16$, we will use
$1+ {h \over \delta\log R} < {2h \over \delta\log R}$. Then, from
(3.17), (3.21) and (3.23), we have
\begin{equation}
\sum_{\scriptstyle N < p_{j} \leq 2N \atop \scriptstyle p_{j+1}-p_j
\leq h}1 \gtrsim  {N\over \log N}{1\over 2^{{4k\log 3N\over
\delta\log R}}} {{2\ell \choose \ell}^2 (k+\ell)!^4 \delta^{2k+2}
\over \eta(k+2\ell)!^2 (2k+2)!^2 2^{k+10}}. \label{eq: 3.24}
\end{equation}
With the values specified in (3.22), the dominating factor in the
coefficient on the right-hand side of (3.24) is
\begin{equation}
2^{-{4k\log 3N\over \delta\log R}} > e^{-65 \bigl({4\over
\eta}\bigr)^6 \log 2}. \label{eq: 3.25}
\end{equation}
The other factors in (3.24) give rise to exponents which are
$\displaystyle O({1\over \eta^{2}}\log{1\over \eta})$. Thus we
obtain \vskip .1in

\noindent {\bf Theorem 2} \,\, {\it For sufficiently small but fixed
$\eta >0$,}
\begin{equation}
\sum_{\scriptstyle N < p_{j} \leq 2N \atop \scriptstyle p_{j+1}-p_j
\leq \eta\log N}1 \gtrsim e^{-c_{3}\eta^{-6}}{N\over \log N}, \qquad
(N\to \infty), \label{eq: 3.26}
\end{equation}
{\it where we can take $c_{3}= \lceil 65 \cdot 4^6 \cdot \log 2
\rceil = 184544$.} \vskip .1in

(This is not the strongest estimate the present method yields.
Taking $\delta = {1\over \ell^{3+c_4}}$ with any fixed $c_{4}>
{1\over 2}$, leads to an estimate of the type (3.26) with
$\eta^{-(5+c_{4})}$ instead of $\eta^{-6}$).

Note that keeping $Q(N,h)$ all the way down to (3.24), by first
keeping it in (3.17) instead of writing the right-hand side of
(3.17) via (3.4), yields $Q(N,h) \gg N$, meaning that the proportion
of natural numbers $n \in [N, 2N]$ for which one can find at least
two primes within a distance of $h$ from $n$ is positive.

\section{Sparsity of very small gaps between primes}

The following result expresses that very small gaps between
consecutive primes occur rarely, in the sense that such gaps do not
constitute a positive proportion of all gaps between consecutive
primes. \vskip .1in

\noindent {\bf Theorem 3} \,\, {\it For any $h > 2$, as
$x\to\infty$, we have}
\begin{equation}
\#\{p_n \leq x \, ; \, p_{n+1}-p_n \leq h \} \ll \min({h\over \log
x},1)\, \pi(x). \label{eq: 4.1}
\end{equation}
{\it In particular, if $h = o(\log x)$, then}
\begin{equation}
\#\{p_n \leq x \, ; \, p_{n+1}-p_n \leq h \} = o(\pi(x)). \label{eq:
4.2}
\end{equation}
\vskip .1in \noindent We remark that for $h=\eta\log x,\, 0< \eta <
1$, the upper estimate for the density of small gaps given by (4.1)
corresponds to the conjectured density $1-e^{-\eta}$ from (1.3)
apart from the constant implied by the $\ll$ symbol; i.e. the simple
upper estimate argument in the proof given below is optimal except
for this constant. \vskip .1in \noindent {\it Proof:} \,\,\, Given
two prime numbers $p,p'$ satisfying $0 < p' - p \leq h$, let us
write $u+h_1 =p,\, u+h_2 = p'$. There are $h$ ordered pairs
$(u,h_1)$ with $1\leq h_1 \leq h$ such that $u+h_1 =p$, and for any
ordered pair $(u,h_1)$ the value of $h_2$ with $h_1 < h_2 \leq 2h$
is fixed. Hence we see that
\begin{align}
h \sum_{\scriptstyle N < p,p' \leq 2N \atop \scriptstyle 0< p'-p
\leq h}1 & < \sum_{\scriptstyle 1\leq h_1 , h_2 \leq 2h \atop
\scriptstyle h_1 \neq h_2}\sum_{\scriptstyle{N\over 2} < u < 3N
\atop \scriptstyle u+h_1 , \, u+h_2 : \,\mathrm{prime}}1 \nonumber \\
& \ll \sum_{\scriptstyle 1 < h_1 , h_2 \leq 2h \atop h_1 \neq h_2}
\mathfrak S(\{h_1,h_2\}){N\over \log^{2}N} \nonumber  \\
& \ll {h^2 N\over \log^{2}N}, \label{eq: 4.3}
\end{align}
where we have used the well-known (see Theorem 5.7 of \cite{HR} or
Theorem  4 in \S 2.3.3 of \cite{Gr}) sieve bound for prime tuples
\begin{equation}
\sum_{N < n \leq 2N}\theta(n+h_1) \cdots \theta(n+h_k) \lesssim 2^k
k! \mathfrak{S}(\mathcal{H})N \label{eq: 4.4}
\end{equation}
with $k=2$, and Gallagher's result (1.16). Thus we have obtained
\begin{equation}
\sum_{\scriptstyle N < p,p' \leq 2N \atop \scriptstyle 0< p'-p \leq
h}1 \ll {hN\over \log^{2}N} \ll {h\over \log N}\pi(N). \label{eq:
4.5}
\end{equation}
We also note that (4.4) used with $k=3$ shows that of the $p,p'$ in
(4.5), the number of those which are not consecutive is $\ll
\bigl({h\over \log N}\bigr)^2 \pi(N)$.

\section{Conditional results}

For the circumstance specified by (3.2) we shall now consider the
consequence of assuming that the level of distribution of primes
$\vartheta$ is greater than ${1\over 2}$. The conditions of
Propositions 1 and 2 allow us to take
\begin{equation}
R= N^{{\vartheta -\epsilon\over 2(1+\delta)}} \label{eq: 5.1}
\end{equation}
with $\epsilon$ and $\delta$ arbitrarily small fixed positive
numbers. We let
\begin{equation}
\ell = \lfloor {\sqrt{k}\over 2}\rfloor . \label{eq:5.2}
\end{equation}
For a given $\vartheta > {1\over 2}$, we determine $k=k(\vartheta)$
sufficiently large and $\epsilon$ and $\delta$ small enough so as to
ensure that the quantity ${k\over k+2\ell +1}{2(2\ell +1)\over \ell
+1}{\log R\over \log N}-1$ occurring in (3.21) is positive. Now $k$
is not necessarily large enough to satisfy (2.8) and the
corresponding inequality when there is the twisting with primes, so
instead of the error term $O(k^3 \delta^2)$ in Propositions 1 and 2,
and in (3.21) we will have the cruder $O_{k}(\delta)$ (or else we
can re-do the calculation as of (2.4) up until (2.8) without having
error terms in what will correspond to (2.5) and (2.6), but this
won't be necessary for our purpose). By choosing a smaller $\delta$
if necessary, we will have the factor in brackets in (3.21) (with
$0$ in place of $\eta$) greater than a positive quantity which
ultimately depends only on $\vartheta$. Hence, comparing (3.17) and
(3.21) we immediately obtain \vskip .1in

\noindent {\bf Theorem 4} \,\, {\it Assume that the primes satisfy a
level of distribution $\vartheta > {1\over 2}$. Let $\eta$ be a
fixed positive small number. Then there exists an integer
$k(\vartheta)$ and a constant $c_{5}(\vartheta)$ such that}
\begin{equation}
\sum_{\scriptstyle N < p_{j} \leq 2N \atop \scriptstyle p_{j+1}-p_j
\leq \eta\log N}1 \gtrsim
c_{5}(\vartheta)\eta^{k(\vartheta)-1}{N\over \log N}, \qquad (N\to
\infty). \label{eq: 5.3}
\end{equation}
\vskip .1in

Notice that the unconditional estimate (3.26) in which $\eta$ takes
place exponentially, gets improved to estimates involving just
powers of $\eta$ when it is assumed that the primes satisfy a level
of distribution greater than ${1\over 2}$. By comparing the factor
in brackets in (3.21) with the corresponding factor in the argument
in \cite{GPY1}, we see that the smallest possible $k(\vartheta)$ we
can assert is either the smallest $r=r(\vartheta)$ such that every
admissible $r$-tuple is guaranteed by the proof of Theorem 1 of
\cite{GPY1} to contain at least two primes infinitely often or it is
$r+1$ (depending on the value of $\vartheta$). The greater the level
of distribution, the smaller power of $\eta$ will be needed in
(5.3). A table of values of $r(\vartheta)$ was provided between
(3.4) and (3.5) of \cite{GPY1} (to avoid confusion we have called
the $k$ in that table as $r$ here). Thus, if $\vartheta
> {20\over 21}$, then we can take $k=7,\, \ell=1$, so that the
$\eta$-dependent factor in right-hand side of (5.3) is $\eta^6$.
However, we recall that assuming $\vartheta \geq 0.971$ and by
considering a linear combination of the $\Lambda_{R}(n;
\mathcal{H},\ell)$ with $k=6$ and $\ell = 0, 1$, the argument for
proving Theorem 1 of \cite{GPY1} still works, so that under this
assumption we can get a lower bound in (5.3) which has $\eta^5$. We
also see from (1.3) that the true order of magnitude of the
$\eta$-dependent factor in right-hand side of (5.3) is believed to
be $\eta$. Thus for this argument to lead to the true order of
magnitude we need to be able to work with admissible pairs
($2$-tuples). But this seems to require improving the results of
\cite{GPY1} to the extent of proving the twin prime hypothesis under
the Elliott-Halberstam conjecture. \vskip .1in

When $\vartheta$ is slightly greater than ${1\over 2}$, from the
condition in Proposition 2, we write
\begin{equation}
R= N^{{{1\over 2}+\xi\over 2(1+\delta)}}, \label{eq: 5.4}
\end{equation}
where we assume that $\xi > 0$ is small.  We take
\begin{equation}
k=2(\ell +1)(2\ell +1), \quad \delta= {1\over \ell^{4}}, \label{eq:
5.5}
\end{equation}
so that
\begin{equation}
\Big\{{k\over k+2\ell+1}{2(2\ell +1)\over \ell +1}{\log R\over \log
N} + \eta -1 +O(k^3 \delta^2)\Big\} = \eta + 2\xi -{1\over \ell} -
{2\xi\over \ell} + O\bigl({1\over \ell^{2}}\bigr). \label{eq: 5.6}
\end{equation}
For a given $\xi$, we determine $\ell$ by
\begin{equation}
\ell = \Big\lceil{1\over \xi}\Big\rceil , \label{eq: 5.7}
\end{equation}
and then the quantity in (5.6) is $>\eta + {\xi\over 2}$ if $\xi$ is
sufficiently small. Hence from (3.21) we now have
\begin{equation}
\tilde{S}_{R} > {{2\ell \choose \ell}\over (k+2\ell)!}(\log N) (\log
R)^{2\ell}(\eta + {\xi\over 2}).\label{eq: 5.8} \end{equation} As
before, we derive an upper-bound for $\tilde{S}_{R}$ starting from
(3.7), together with (3.8), (3.13) and (3.15). In our case $u={h
\over \delta\log R}$, and upon using the relations in (3.2), (5.4),
(5.5) and (5.7), we have
\begin{equation}
u= {4\Bigl(1+ \Big\lceil{1\over \xi}\Big\rceil^4 \Bigr)\over
1+2\xi}\eta. \label{5.9} \end{equation} This is a small quantity if
for a given small $\xi$ we take $\eta$ small enough, say $\eta \leq
{\xi^4\over 5}$, so that we can say $(1+u)^{k+2} < 2^{k+2}$. Using
this in (3.15) and (3.13) gives
\begin{equation}
I < N (\log R)^{4(k+\ell)}(\log N)^2 \Bigl({h \over \delta\log
R}\Bigr)^k \,\, {2^{{4k\log 3N\over \delta\log R}+k+2}(2k+2)!^2\over
(k+\ell)!^4}. \label{eq: 5.10}
\end{equation}
Plugging this in (3.7), and using that together with (3.4) and (5.8)
we obtain
\begin{equation}
\sum_{\scriptstyle N < p_{j} \leq 2N \atop \scriptstyle p_{j+1}-p_j
\leq h}1 \gtrsim  {N\over \log N}\eta^{k-1}\Bigl({\xi\over
3}\Bigr)^2 \Bigl({4\delta(1+\delta)\over 1+2\xi}\Bigr)^k {{2\ell
\choose \ell}^2 (k+\ell)!^4 k!^2 \over (k+2\ell)!^2 (2k+2)!^2
2^{{4k\log 3N\over \delta\log R}+k+2}}. \label{eq: 5.11}
\end{equation}
By the relations (5.4), (5.5) and (5.7), all of the factors after
$\eta^{k-1}$ can be expressed in terms of $\xi$. What interests us
most is the power of $\eta$, so we re-express (5.11) as
\begin{equation}
\sum_{\scriptstyle N < p_{j} \leq 2N \atop \scriptstyle p_{j+1}-p_j
\leq \eta\log N}1 \gtrsim
c_{6}(\xi)\eta^{4\xi^{-2}+14\xi^{-1}+11}{N\over \log N}, \qquad
(N\to \infty), \label{eq: 5.12}
\end{equation}
valid when the level of distribution of primes is assumed to allow
us to take $R$ as in (5.4) which can be re-written as
\begin{equation}
R = N^{{1+2\xi\over 4(1+\lceil {1\over \xi}\rceil^{-4})}},
\label{5.13}
\end{equation}
for fixed $\eta \in (0,{\xi^4\over 5}]$. Here $\xi$ has to be
sufficiently small, which ensures that $\ell$ and $k$ are
sufficiently large so as to permit the inequality (5.8). \vskip .1in

In \S 3 of \cite{GPY1} it was shown that under the
Elliott-Halberstam conjecture we have
\begin{equation}
\liminf_{n\to\infty}{p_{n+2}-p_{n}\over \log p_{n}} =0. \label{eq:
5.14}
\end{equation}
Our method also shows that such gaps occur in positive proportion.
To see this we consider
\begin{equation}
\tilde{S}_{R,2} := {1\over N(h\log R)^{k}}\sum_{N < n \leq
2N}\bigl(\Theta(n,h) - 2\log 3N\bigr) \bigl(\sumstar_{\mathcal H}
\Lambda_{R}^{2}(n;\mathcal H , \ell)\;\bigr). \label{eq: 5.15}
\end{equation}
As was done in \cite{GPY1} along with the modification provided by
(2.8), we find
\begin{equation}
\tilde{S}_{R,2} \gtrsim {{2\ell \choose \ell}\over (k+2\ell)!}(\log
N) (\log R)^{2\ell}\Big\{{k\over k+2\ell+1}{2(2\ell +1)\over \ell
+1}{\log R\over \log N} + \eta -2 - k^3 \delta^2 \Big\}. \label{eq:
5.16}
\end{equation}
Here we are assuming that $\vartheta =1$, and so we can take $R=
N^{{1\over 2(1+\delta)}}$. In the proof of (5.14), $k$ is taken to
be sufficiently large, $\ell =\lfloor {\sqrt{k}\over 2} \rfloor$. If
$\delta$ is taken to be accordingly small, say $\delta = {1\over
k^2}$, then the quantity in brackets in (5.16) is
\begin{align}
& > 2\bigl(1-{2\ell +1\over k}-{1\over 2\ell}\bigr)(1-\delta)+\eta
-2 -k^{2}\delta^{3} \nonumber \\
& > \eta - {2(2\ell +1)\over k}-{1\over \ell} -2\delta -k^{2}\delta^{3} \nonumber \\
& > \eta -{2(\sqrt{k}+1)\over k}-{2\over \sqrt{k}-2}-{2\over
k^2}-{1\over k} \nonumber \\
& > \eta -{5\over \sqrt{k}} -{3\over k}-{2\over k^2}  \qquad ({\rm
for} \,\,\, k>36) \nonumber \\
& > \eta - {6\over \sqrt{k}} \nonumber \\
& > {\eta\over 2} \qquad ({\rm for} \,\,\, k>{144\over \eta^2}).
\label{eq:5.17}
\end{align}
The rest of the argument is almost identical to what was done as of
(3.7), the only changes are that we now have the summation condition
$\Theta(n,h) \geq {5\over 2}\log N$, and ${h\over \delta\log R} =
2\eta\bigl(1+{1\over \delta}\bigr)$ being not small we should use
some bound like $(1+u)^{k+2} \leq (2u)^{k+2}$ (cf. between (5.9) and
(5.10)). The following is the result of this calculation. \vskip
.1in

\noindent {\bf Theorem 5} \,\, {\it Assuming the Elliott-Halberstam
conjecture we have}
\begin{equation}
\sum_{\scriptstyle N < p_{j} \leq 2N \atop \scriptstyle p_{j+2}-p_j
\leq \eta\log N}1 \gtrsim e^{c_{7}\eta^{-2}\log\eta}{N\over \log N},
\qquad (N\to \infty) \label{eq: 5.18}
\end{equation}
($c_7 = 5$ gives a valid result if $\eta$ is small enough).

\newpage
\footnotesize D. A. Goldston  \,\,\,
(goldston@math.sjsu.edu)

Department of Mathematics

San Jose State University

San Jose, CA 95192

 USA \\

J. Pintz \,\,\, (pintz@renyi.hu)

R\'enyi Mathematical Institute of the Hungarian Academy of Sciences

H-1364 Budapest

P.O.B. 127

Hungary \\

C. Y. Y{\i}ld{\i}r{\i}m \,\,\, (yalciny@boun.edu.tr)

Department of Mathematics

Bo\~{g}azi\c{c}i University

Bebek, Istanbul 34342

Turkey \\


\vspace{1cm}

\begin{thebibliography}{20}
\bibitem{BD}
E. Bombieri and H. Davenport, Small differences between prime
numbers, {\it Proc. Roy. Soc. Ser. A} {\bf 293} (1966), 1-18.

\bibitem{Ga}
P. X. Gallagher, On the distribution of primes in short intervals,
{\it Mathematika} {\bf 23} (1976), 4--9.

\bibitem{GPY1}
D. A. Goldston, J. Pintz and C. Y. Y{\i}ld{\i}r{\i}m, Primes in
tuples~I, {\it Ann. of Math.} {\bf 170} (2009), 819--862.

\bibitem{GPY2}
D. A. Goldston, J. Pintz and C. Y. Y{\i}ld{\i}r{\i}m, Primes in
tuples~II, {\it Acta Math.} {\bf 204} (2010), 1--47.

\bibitem{GPY}
D. A. Goldston, J. Pintz and C. Y. Y{\i}ld{\i}r{\i}m, Positive
proportion of small gaps between consecutive primes, {\it Publ.
Math. Debrecen} (to appear).

\bibitem{GY1}
D. A. Goldston and C. Y. Y{\i}ld{\i}r{\i}m, Higher correlations of
divisor sums related to primes I: Triple correlations, {\it
Integers} {\bf 3} (2003), A5, 66pp. (electronic).

\bibitem{GY}
D. A. Goldston and C. Y. Y{\i}ld{\i}r{\i}m, Higher correlations of
divisor sums related to primes III: Small gaps between primes, {\it
Proc. London Math. Soc. (3)} {\bf 95} (2007), no. 3, 653-686.

\bibitem{Gr} G. Greaves, {\it Sieves in number theory}, Springer,
Berlin, 2001.

\bibitem{HR} H. Halberstam and H.-E. Richert, {\it Sieve methods},
Academic Press, New York, 1974.

\bibitem{H}
M. N. Huxley, Small differences between consecutive primes II, {\it
Mathematika} {\bf 24} (1977), 142-152.

\bibitem{M}
H. Maier, Small differences between prime numbers, {\it Michigan
Math. J.} {\bf 35} (1988), 323-344.

\bibitem{P}
J. Pintz, Are there arbitrarily long arithmetic progressions of twin
primes?, in {\it An irregular mind, Szemer\'{e}di is 70}, (Editors:
I. Barany and J. Solymosi), 525-559, Bolyai Soc. Math. Studies {\bf
21}, Springer, 2010.

\bibitem{RS}
J. B. Rosser and L. Schoenfeld, Approximate formulas for some
functions of prime numbers, {\it Illinois J. Math.} {\bf 6} (1962),
64-94.

\bibitem{So1}
K. Soundararajan, Small gaps between prime numbers: The work of
Goldston-Pintz-Y{\i}ld{\i}r{\i}m, {\it Bull. Amer. Math. Soc.
(N.S.)} {\bf 44} (2007), 1-18.

\bibitem{So2}
K. Soundararajan, The distribution of prime numbers. {\it
Equidistribution in number theory, an introduction}, 59--83, NATO
Sci. Ser. II Math. Phys. Chem., 237, Springer, Dordrecht, 2007.

\end{thebibliography}
\end{document}